\documentclass[a4paper,12pt]{article}
\usepackage{etex}

\usepackage[utf8x]{inputenc}
\usepackage{arabtex}
\usepackage{afterpage}
\usepackage{hyperref}
\usepackage{relsize}
\usepackage{mathrsfs}
\usepackage[all]{xy}

\usepackage{tikz}
\usetikzlibrary{matrix}

\usepackage{amsmath,amssymb,amsfonts,amscd, graphicx, latexsym, verbatim, multirow, color}
\usepackage{epstopdf}
\usepackage{epsfig}
\usepackage{geometry} % to change the page dimensions
\usepackage{adjustbox}

\newcommand{\Q}{\mathbb Q}
\newcommand{\Z}{\mathbb Z}

 \renewcommand{\R}{\mathbf R}

\newcommand{\pgl}{\mathrm{PGL}_2(\Z)}

\usepackage{color, ulem}

\newcommand{\Jimm}{\mathbf J}
\newcommand{\jimm}{{\adjustbox{scale=.5, raise=0.6ex, trim=0px 0px 0px 7px, padding=0ex 0ex 0ex 0ex}{\RL{j}}}}

\newcommand{\sherhhh}[1]{\fboxsep=0pt\setlength{\fboxrule}{1pt}
\begin{center}
   \fbox{\colorbox{red}{
         \begin{minipage}[t]{13cm}
            #1
         \end{minipage}
      }
   }
\end{center}}

\renewcommand{\sherhhh}[1]{}
\begin{document}

%opening
\title{Dynamics of a family of continued fraction maps}
\author{Muhammed Uluda\u{g}$^*$, Hakan Ayral\footnote{
{Department of Mathematics, Galatasaray University}
{\c{C}{\i}ra\u{g}an Cad. No. 36, 34357 Be\c{s}ikta\c{s}}
{\.{I}stanbul, Turkey}}}

\maketitle

%\pagecolor{yellow}

\begin{abstract}
We study the dynamics of a family of continued fraction maps parametrized by the unit interval. 
This family contains as special instances the Gauss continued fraction map and the Fibonacci map.
We determine the transfer operators of these dynamical maps and make a preliminary study of them.
We show that their analytic invariant measures obeys a common functional equation generalizing Lewis' functional equation and we find invariant measures for some members of the family.   
We also discuss a certain involution of this family which sends the Gauss map to the Fibonacci map.
 \end{abstract}
\section{Introduction}\label{sec:introduction}
Denote the continued fractions in the usual way
$$
[0,n_1,n_2,\dots]=\cfrac{1}{n_1+\cfrac{1}{n_2+\cfrac{1}{\dots}}},
$$
where $n_i{> 0}$ $(i=1,2\dots)$ are integers.
Our aim in this paper is to study a family of continued fraction maps ${\mathbb T}_\alpha$, parametrized by $\alpha\in [0,1]$ and defined as
\begin{equation}\label{defofalpha}
{\mathbb T}_\alpha(x)=
\begin{cases} 
[0,m_{k+1}, m_{k+2},m_{k+3},\dots]&n_k>m_k \,(*)\\
[0,m_{k}-n_{k}, m_{k+1},m_{k+2},\dots]&n_k<m_k\, (**)
\end{cases}
\end{equation}
where
$\alpha=[0,n_1,n_2, \dots]$ and $x=[0,m_1, m_2, \dots].$

\begin{center}
\noindent{\includegraphics[scale=.50]{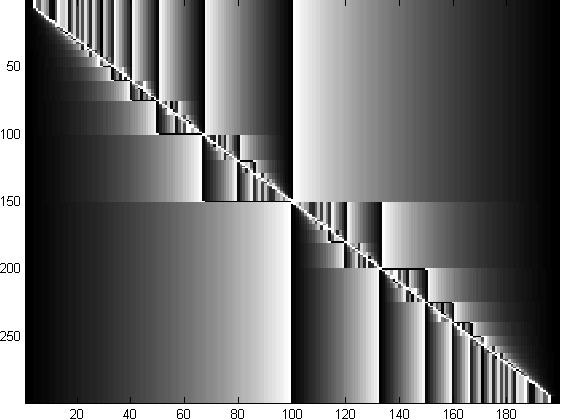}}\\
{\small Plot of the  map $ {\mathbb T}_\alpha(x)$ as a function of $\alpha$ and $x$. 
The intensity is proportional to the value of $T_\alpha(x)$. The symmetry reflects the fact that 
$T_{1-\alpha}(1-x)=T_\alpha(x)$.}
\end{center}
The first thing to notice about this family is that it is not continuous in the variable $\alpha$,  the two representations of rational $\alpha$ as a continued fraction leads to two different maps with distinct dynamical behaviors.
Another problem is that ${\mathbb T}_\alpha(x)$ is not defined if $x=\alpha$ or if $x\in \mathbb Q$ with one of its two continued fraction expansions being an initial segment of $\alpha$. We overcome this issue by assuming that the rational $x$ are represented by continued fractions $[0,m_1, m_2, \dots, m_k, \infty]$ with $m_k>1$. 
Consequently, we set ${\mathbb T}_\alpha(x)=0$ for such $x$.
We also set ${\mathbb T}_\alpha(\alpha):=0$.

Each ${\mathbb T}_\alpha(x)$ constitutes a M\"obius system \cite{schweigerpaper} which is a special kind of a fibered system \cite{schweigerbook}. For $\alpha=[0,\infty]$, the map ${\mathbb T}_\alpha$ is the Gauss continued fraction map \cite{kramp}, \cite{karma} and for 
$\alpha=[0,1,1,\dots]$, ${\mathbb T}_\alpha$ is the so-called Fibonacci map \cite{isola2014continued}. In Section 3 below we recall some basic facts about them and prove that they are in fact conjugates of each other. In the preceding Section 2, we define 
a transfer operator for ${\mathbb T}_\alpha$ and prove that their
invariant functions obeys a common functional equation generalizing Lewis's. Section 4 is devoted to the study of some other special instances of the family ${\mathbb T}_\alpha$. In particular for $\alpha=[0,K,K,\dots]$ ($K=1,2,\dots$) we determine an analytic invariant measure for  ${\mathbb T}_\alpha$.

\section{Transfer operators and functional equations}\label{family}
To give another perspective on  $\mathbb T_\alpha$, 
first represent  $\alpha=[0,n_1,n_2,\dots]$ by a binary string
$\mathbf 0_{n_1}\mathbf 1_{n_2} \mathbf 0_{n_3}\dots $ (always starting with a $\mathbf 0$) with $n_i>0$ for all $i$.
Also represent $x=[0,m_1, m_2, \dots]$ in this form. Then $\mathbb T_\alpha$ compares the strings of $x$ and $\alpha$ and returns $x$ with 
the initial segment matching with that of $\alpha$ chopped off. If the remaining sequence starts with a $1$, then it
permutes $0$'s with $1$'s.\footnote{Alternatively, one may assume that $\mathbf 0_{n_1}\mathbf 1_{n_2} \mathbf 0_{n_3}\dots $ and $\mathbf 1_{n_1}\mathbf 0_{n_2} \mathbf 1_{n_3}\dots $ represents the same continued fraction.}

\subsection{The transfer operator}
To write a transfer operator for  ${\mathbb T}_\alpha$, note that its inverse branches are
$$
{\mathbb T}_\alpha^{-1}(y)=
\begin{cases} [0, n_1, n_2, \dots, n_{k-1}, i+y] &(1\leq k; 1\leq i<n_{k}), \mbox{ or }\\
 [0, n_1, n_2, \dots, n_{k-1}, n_k,y] &(1\leq k).
 \end{cases}
$$
Set $C_\alpha(\psi):=\psi \circ {\mathbb T}_\alpha$ for the Koopman operator. 
Its left-adjoint with respect to its action on the set of cumulative distribution functions of  probability measures on $[0,1]$ is the ``{integral}" Gauss-Kuzmin-Wirsing operator
\begin{eqnarray*}
({\mathbb L}_{s, \alpha} \Psi)(y)=
\sum_{k=1}^\infty (-1)^{1+k}
\bigl(\Psi[0, n_1, \dots, n_{k-1}, n_k,y]-\Psi[0, n_1, \dots, n_{k-1}, n_k,0]\bigr)\qquad
\\
+\sum_{k=1}^\infty (-1)^{k}
\sum_{i=1}^{n_k-1} \bigl(\Psi[0,n_1,\dots, n_{k-1},i+y]-\Psi[0,n_1,\dots, n_{k-1},i]\bigr)
\end{eqnarray*}
Assuming $\Psi$ is differentiable and the convergence is nice, we may take derivatives which yields the operator
\begin{eqnarray*}
{\mathscr L}_{s=1, \alpha} \psi=
\sum_{k=1}^\infty (-1)^{1+k}\left\{\frac{\mathrm d}{{\mathrm d}y}
[0,n_1,\dots, n_{k-1},n_k,y]\right\}\psi[0, n_1, \dots, n_{k-1}, n_k,y]+\qquad\\
\sum_{k=1}^\infty (-1)^{k}
\sum_{i=1}^{n_k-1} \left\{\frac{\mathrm d}{{\mathrm d}y}
[0,n_1,\dots, n_{k-1},i+y]\right\} \psi[0,n_1,\dots, n_{k-1},i+y]
\end{eqnarray*}
We will call this the {\it Gauss-Kuzmin-Wirsing (GKW) operator} of ${\mathbb T}_\alpha$.

To handle the alternation of  signs, set 
$$
[0, n_1, \dots, n_{k-1}, n_k,y]=:\frac{A_k y+ B_k}{C_k y+D_k},
$$
where  $A_k, B_k, C_k, D_k\in \Z$ and $A_kD_k-B_kC_k=(-1)^{1+k}$. Then 
$$
\frac{\mathrm d}{{\mathrm d}y}
[0,n_1,\dots, n_{k-1},n_k,y]=\frac{(-1)^{k+1}}{(C_k y+D_k)^2}.
$$
Now we define the transfer operator as below, which specializes to the Gauss-Kuzmin-Wirsing operator when  $s=1$:
\begin{eqnarray*}
{\mathscr L}_{s, \alpha} \psi:=
\sum_{k=1}^\infty \frac{1}{(C_k y+D_k)^{2s}} \psi[0, n_1, \dots, n_{k-1}, n_k,y]+\qquad\qquad\qquad\qquad\qquad\qquad\\
\qquad\qquad\qquad\qquad
\sum_{k=1}^\infty \sum_{i=1}^{n_k-1} \frac{1}{(C_{k-1} (y+i)+D_{k-1})^{2s}} \psi[0,n_1,\dots, n_{k-1},i+y]
\end{eqnarray*}
We will consider $y$ as a real variable and use the following practical form for ${\mathscr L}_{s, \alpha}$:
\begin{eqnarray*}
{\mathscr L}_{s, \alpha} \psi=
\sum_{k=1}^\infty \left|\frac{\mathrm d}{{\mathrm d}y}
[0,n_1,\dots, n_{k-1},n_k,y]\right|^{s} \psi[0, n_1, \dots, n_{k-1}, n_k,y]+\qquad\qquad\\
\sum_{k=1}^\infty 
\sum_{i=1}^{n_k-1} \left|\frac{\mathrm d}{{\mathrm d}y}
[0,n_1,\dots, n_{k-1},i+y]\right|^{s} \psi[0,n_1,\dots, n_{k-1},i+y]
\end{eqnarray*}
The first summand can be written as
\begin{eqnarray}\label{rationalalpha}
\sum_{k=1}^\infty \left|\frac{\mathrm d}{{\mathrm d}y}
[0,n_1,\dots, n_{k-1},y]\right|^{s} \psi[0, n_1, \dots, n_{k-1},y]
-
 \left|\frac{\mathrm d}{{\mathrm d}y}
[0,y]\right|^{s} \psi[0,y],
\end{eqnarray}
where for rational $\alpha=[0,n_1,\dots, n_i]$, the terms $k>i$ in (\ref{rationalalpha}) vanish.\footnote{We stress that for rational $\alpha$ the two different choices for the continued fraction representation of $\alpha$ 
will lead to two different dynamical maps and transfer operators. We will discuss these below.}
Hence $({\mathscr L}_{s, \alpha} \psi)(y)=$
\begin{eqnarray}\label{domain}
-\frac{1}{y^{2s}}\psi\left(\frac{1}{y}\right)+
\sum_{k=1}^\infty 
\sum_{i=0}^{n_k-1} \left|\frac{\mathrm d}{{\mathrm d}y}
[0,n_1,\dots, n_{k-1},i+y]\right|^{s} \psi[0,n_1,\dots, n_{k-1},i+y]
\end{eqnarray}
A natural domain for ${\mathscr L}_{s, \alpha} $ is the space $BV[0,1]$ of functions of bounded variation on the interval $[0,1]$. The first term in (\ref{domain}) apparently violates this and requires that $\psi$ to be defined on a larger domain. 
However, this term serves only to counterbalance its identical twin hidden in the infinite summand in (\ref{domain}).

Consider the action of $\pgl$ on the set of functions on $\R$ defined by 
$$
M(x)=\frac{ax+b}{cx+d}\in\pgl\implies (\psi|M)=\frac{1}{|cx+d|^{2s}}\psi(\frac{ax+b}{cx+d})
$$
Set $T(x):=x+1$ and $U(x):=1/x$.
Then
\begin{eqnarray*}
 \left|\frac{\mathrm d}{{\mathrm d}y}
[0,n_1,\dots, n_{k-1},i+y]\right|^{s} \psi[0,n_1,\dots, n_{k-1},i+y]=(\psi|UT^{n_1}U\dots UT^{n_{k-1}}UT^{i})(y)
\end{eqnarray*}
and one may write
\begin{eqnarray*}
{\mathscr L}_{s, \alpha} \psi=
\sum_{k=1}^\infty (\psi|UT^{n_1}U\dots UT^{n_{k}}U)+
\sum_{k=1}^\infty 
\sum_{i=1}^{n_k-1} 
(\psi|UT^{n_1}U\dots UT^{n_{k-1}}UT^{i}).
\end{eqnarray*}
Alternatively,
\begin{eqnarray}
{\mathscr L}_{s, \alpha} \psi=
-(\psi|U)+
\sum_{k=1}^\infty 
\sum_{i=0}^{n_k-1} (\psi|UT^{n_1}U\dots UT^{n_{k-1}}UT^{i})\label{first}\\
=
(\psi|-U+
\sum_{k=1}^\infty 
\sum_{i=0}^{n_k-1} UT^{n_1}U\dots UT^{n_{k-1}}UT^{i})\\
=
(\psi|-U+
\sum_{k=1}^\infty 
UT^{n_1}U\dots UT^{n_{k-1}}U\sum_{i=0}^{n_k-1} T^{i})\\
=
(\psi|-U+
\sum_{k=1}^\infty 
UT^{n_1}U\dots UT^{n_{k-1}}U(I-T)^{-1}(I-T^{n_k})).\label{last}
\end{eqnarray}
Now one has, using the form (\ref{first}) above,
\begin{eqnarray*}
({\mathscr L}_{s, \alpha} \psi|T)=
-(\psi|UT)+
\sum_{k=1}^\infty 
\sum_{i=0}^{n_k-1} (\psi|UT^{n_1}U\dots UT^{n_{k-1}}UT^{i+1})
\end{eqnarray*}
Hence $({\mathscr L}_{s, \alpha} \psi |I-T)$=
\begin{eqnarray*}
-(\psi|U(I-T))+
\sum_{k=1}^\infty 
\bigl[ 
(\psi|UT^{n_1}U\dots UT^{n_{k-1}}U)- (\psi|UT^{n_1}U\dots UT^{n_{k-1}}UT^{n_k})
\bigr]
\end{eqnarray*}
and therefore $(({\mathscr L}_{s, \alpha} \psi|I-T)|U)=$
\begin{eqnarray*}
-(\psi|U(I-T)U)+
\sum_{k=1}^\infty 
\bigl[ 
(\psi|UT^{n_1}U\dots UT^{n_{k-1}})- (\psi|UT^{n_1}U\dots UT^{n_{k}}UT^{n_k}U)
\bigr]
\end{eqnarray*}
Summing the last two equations we get $({\mathscr L}_{s, \alpha} \psi|(I-T)(I+U))=$
\begin{eqnarray*}
-(\psi|U-UT+I-UTU)+(\psi|I+U)=(\psi|UT+UTU)
\end{eqnarray*}
Now suppose that $\psi$ is an eigenfunction with eigenvalue $\lambda$, i.e. ${\mathscr L}_{s, \alpha} \psi=\lambda\psi$. Then
\begin{eqnarray}\label{feq}
(\lambda\psi|(I-T)(I+U))=(\psi|UTU+UT)\implies \\
(\psi|(I-T)(I+U))=\frac{1}{\lambda}(\psi|UTU+UT)
\end{eqnarray}
Observe that $(\psi|I-T-TU)=0$ is precisely Lewis' three-term functional equation and $(\star|UTU+UT)$ is Isola's transfer operator of the Farey map.
Written explicitly, this equation looks as follows:
\begin{eqnarray}\label{corporate}
\psi(y)-
\psi(1+y)+
\frac{1}{y^{2s}}\left\{ 
\psi\left(\frac{1}{y}\right)-\psi\left(1+\frac{1}{y}\right)\right\}=\qquad\qquad\qquad\qquad\qquad\\
\nonumber\qquad\qquad\qquad\qquad\qquad\frac{1}{\lambda(1+y)^{2s}}\left\{
\psi\left(\frac{y}{1+y}\right)+\psi\left(\frac{1}{1+y}\right)\right\}
\end{eqnarray}
It is readily seen in the explicit form that any 1-periodic odd function solves Equation \ref{corporate}. The set of all solutions forms a vector space over $\R$. Note that, for a given solution $\psi$ of  Equation \ref{corporate},
the infinite series defining ${\mathscr L}_{s, \alpha}\psi$ may or may not converge; even if it does, $\psi$ may or may not be an eigenfunction of ${\mathscr L}_{s, \alpha}$. See \cite{lewis2001period} for a study of the case $\alpha=0$. 

For $\lambda=\pm1$, one may regroup  (\ref{corporate}) in several interesting ways. For example, 
\begin{eqnarray}\label{feq2}
 (\psi|I-T-\frac{1}{\lambda}UT)+(\psi|I-T-\frac{1}{\lambda}UT|U)=0\\
\iff (\psi|(I-T-\frac{1}{\lambda}UT)(I+U))=0
\end{eqnarray}
Setting 
$$
\phi:=(\psi|I-T-\frac{1}{\lambda}UT)\iff \phi(y)=\psi(y)-\psi(1+y)-\frac{1}{\lambda (1+y)^{2s}}\psi\left(\frac{1}{1+y}\right)
$$ we arrive at the cocycle relation
\begin{eqnarray*}
\implies \phi+(\phi|U)=0 \iff (\phi|I+U)=0
\end{eqnarray*}
In order to solve the equation $\phi=(\psi|I-(I+\frac{1}{\lambda}U)T)$ for $\psi$, we need an inverse for the operator 
$(\star|I-(I+\frac{1}{\lambda}U)T)$. To be more explicit, this is the operator
$$
{\mathscr B}_\lambda: \psi(y)\to \psi(y)-\psi(1+y)-\frac{1}{\lambda}\frac{1}{(1+y)^2}\psi\left(\frac{1}{1+y}\right)
$$
Its kernel equation is precisely the functional equation for analytic eigenfunctions of Mayer's transfer operator; and specialize to the following three-term functional equation when $\lambda=1$. 
\begin{eqnarray}\label{lws}
{\mathscr B}_1\psi= 0 \iff \psi(y)-\psi(1+y)-\frac{1}{(1+y)^2}\psi\left(\frac{1}{1+y}\right)\equiv 0.
\end{eqnarray}
This functional equation is connected to the one studied by Lewis and Zagier \cite{lewis2001period}.
Modulo the above remark about its kernel, one can formally invert the equation $\phi=(\psi|I-(I+\frac{1}{\lambda}U)T)$  and write
\begin{eqnarray*}
\psi=(\phi|(I-(I+\frac{1}{\lambda}U)T)^{-1})=(\phi|\sum_{k=0}^\infty ((I+\frac{1}{\lambda}U)T)^{k})
=\sum_{k=0}^\infty (\phi|( (I+\frac{1}{\lambda}U)T)^{k})
\end{eqnarray*}
Now if $\lambda=1$ then since $(\phi|I+U)=0$, the above inversion formula becomes simply $\psi=\phi$. Taking into account the kernel of $\mathscr B$, we may write 
$$\psi=\phi+\eta,$$
where $\eta$ is a solution of (\ref{lws}) and $\phi$ satisfies $(\phi|I+U)=0$. Note that
$$
\psi=(f|I-U) \mbox{ for some } f \implies (\psi|I+U)=0.
$$

Another rearrangement of Equation \ref{corporate}, which makes explicit the connection with Lewis' functional equation, is
\begin{eqnarray}\label{feq2}
 (\psi|I-T-\frac{1}{\lambda}UTU)+(\psi|I-T-\frac{1}{\lambda}UTU|U)=0\\
\iff (\psi|(I-T-\frac{1}{\lambda}UTU)(I+U))=0
\end{eqnarray}
Setting 
$
\phi^*:=(\psi|I-T-\frac{1}{\lambda}UTU)
$ gives again the cocycle relation $\phi^*+(\phi^*|U)=0 \iff (\phi^*|I+U)=0.$
The equation  $\phi^*=(\psi|I-(T+UTU))=0$ (i.e. the case $\lambda=1$) is precisely the following three-term functional equation  studied by Lewis and Zagier \cite{lewis2001period}:
\begin{eqnarray}\label{lws}
\psi(y)-\psi(1+y)-\frac{1}{(1+y)^{2s}}\psi\left(\frac{y}{1+y}\right)\equiv 0.
\end{eqnarray}

\subsection{Common invariant measure}
The maps ${\mathbb T}_\alpha$ share a common invariant measure, namely Minkowski's measure, whose cumulative distribution function is given as
\begin{eqnarray}
\mbox{\bf ?}([0,n_1,n_2,\dots,n_{k-1},n_k])=\sum_{k=1}^\infty (-1)^{1+k}2^{-n_1-n_2\dots-n_k}.
\end{eqnarray}
Actually, \mbox{\bf ?} is the common invariant measure of a much wider class of maps $T:[0,1]\mapsto [0,1]$ whose inverse branches are  all $\pgl$ on $[0,1]$.
To see the this, suppose that 
the inverse branches of $T$ are $\{\varphi_\beta\}_{\beta=1,2,\dots}$. Then each $\varphi_\beta$ can be written as 
$$
\varphi_\beta(y)=[0,n_1,n_2,\dots,n_{k-1},i+y],
$$
where $0<k, n_1, n_2, \dots$ and $0\leq i$ depends on $\beta$. 
Suppose $X$ is a random variable on the unit interval with  law $\mbox{\bf ?}$ and set $Y:={T}(X)$.
The law $\mathbf{F}_Y$ of $Y$ is
\begin{eqnarray*}
\mathbf{F}_Y(y)=\mathrm{Prob}\{Y\leq y\}=
\mathrm{Prob}\{{T}(X)\leq y\}
=\sum_{\beta}
\bigl|
\mbox{\bf ?}\left(\varphi_\beta(y)\right)-
\mbox{\bf ?}\left(\varphi_\beta(0)\right)
\bigr|\\
=
\sum_{\beta}
\bigl|
\mbox{\bf ?}[0,n_1,n_2,\dots,n_{k-1},i+y]-
\mbox{\bf ?}[0,n_1,n_2,\dots,n_{k-1},i]
\bigr|\\
=\sum_{\beta}\mbox{\bf ?}(y)2^{-(n_1+\dots+n_{k-1}+i)}\implies\\
\mathbf{F}_Y(y)=\mbox{\bf ?}(y)\sum_{\beta}2^{-(n_1+\dots+n_{k-1}+i)},\\
\end{eqnarray*}
and the series of the last line {\it must} sum up to 1, because $\mathbf{F}_Y$ and $\mbox{\bf ?}(y)$ are both probability laws.

\section{Gauss and Fibonacci maps}
In this section, we will discuss the two special instances of the family ${\mathbb T}_\alpha$, namely the Gauss map ${\mathbb T}_0$ 
and the Fibonacci map ${\mathbb T}_{{\Phi^*}}$, where $\Phi^*=[0,1,1,\dots]$.

%%%%%%%%%
\subsection{The Gauss map}
If $\alpha=0=[0,\infty]$  with $n_1=\infty$, then in the definition (\ref{defofalpha}) of 
${\mathbb T}_\alpha$ we are always in the first case $(*)$,  and so ${\mathbb T}_0$ is the well-known {\it Gauss continued fraction map}
\begin{equation}
{\mathbb T}_0(x)=\frac{1}{x}-\left\{\frac{1}{x}\right\}=[0,m_{2}, m_{3},m_{4},\dots],
\end{equation}
where $x=[0,m_{1}, m_{2},m_{3},\dots]$ and $\{ y\}$ denotes the fractional part of $y$. 
One has
$$
({\mathscr L}_{s,\alpha=1}\psi)(y)=\sum_{n=1}^\infty \frac{1}{(y+n)^{2s}}\psi\left(\frac{1}{y+n}\right)
$$
Its invariant functions are solutions of the  functional equation
\begin{equation}\label{itfe}
p(y)=\sum_{n=1}^\infty \frac{1}{(y+n)^{2s}}p\left(\frac{1}{y+n}\right)
\end{equation}
Substituting $y+1$ in place of $y$ leads to the three-term functional equation
\begin{eqnarray}\label{gfe}
p(y)=p(y+1)+\frac{1}{(y+1)^{2s}}p\left(\frac{1}{y+1}\right),
\end{eqnarray}
which admits, when $s=1$ the function (called the {\it Gauss density}) 
$$
p_0(y):={1\over \log 2} {1 \over 1+y}
$$ 
as a solution. 
It is well known and easily verified that the Gauss density is an invariant density for the classical GKW-operator 
${\mathscr L}_{s=1,\alpha=1}$. 

\begin{figure}[h]
\begin{center}
\noindent{\includegraphics[scale=.30]{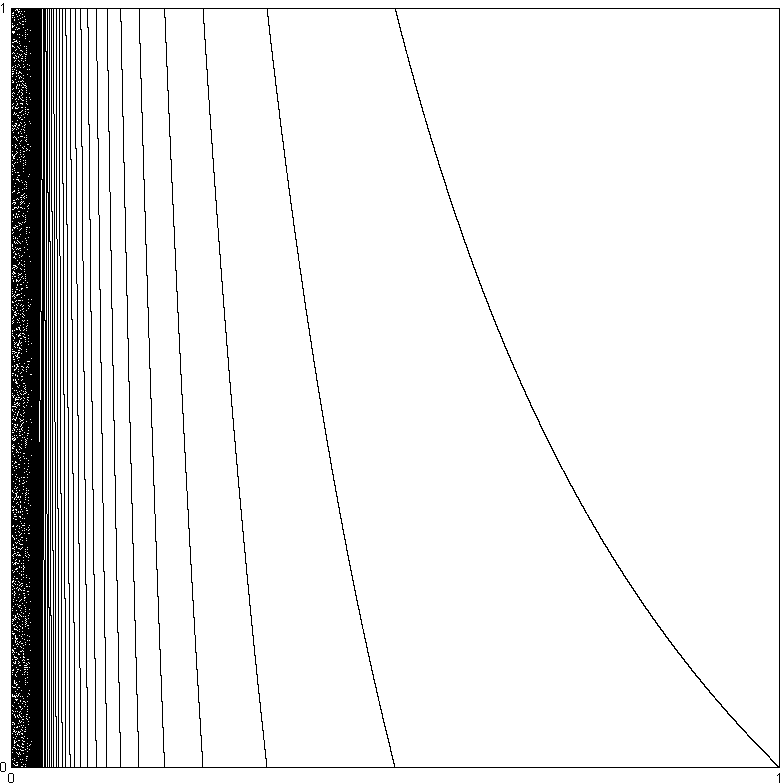}}
\noindent{\includegraphics[scale=.30]{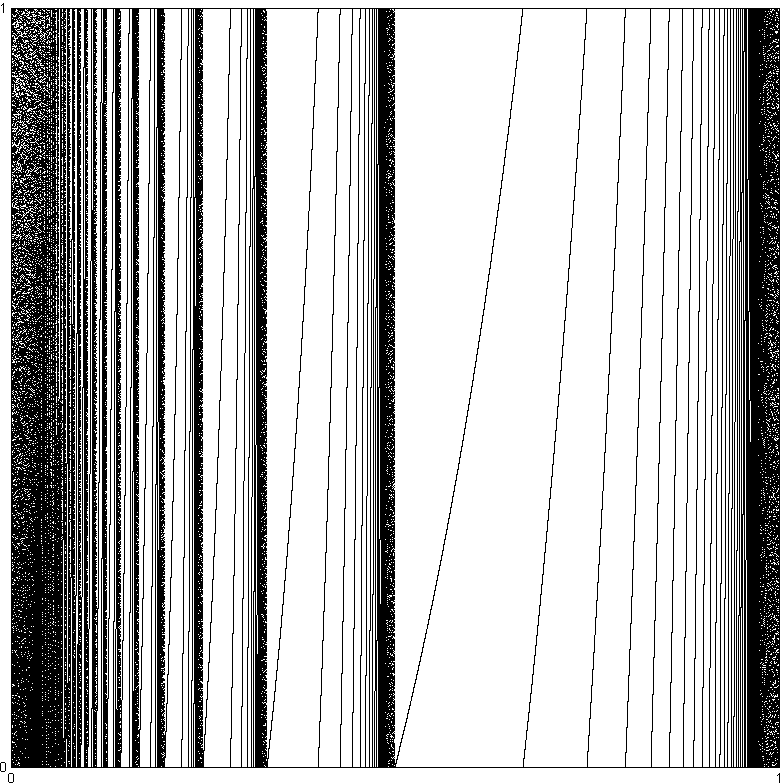}}\\
The Gauss map ${\mathbb T}_0$ and its second iteration.
\end{center}
\end{figure}

Mayer \cite{mayer1990thermodynamic} considered the operator ${\mathscr L}_{s,\alpha=1}$ as acting on the space of functions holomorphic on the disc $|z-1|\leq 3/2$ and proved that the determinant of $1-{\mathscr L}_{s,\alpha=1}$ exists in the Fredholm sense and equals to the Selberg zeta function of the extended modular group.

\subsection{The Fibonacci map}

If  $\alpha={\Phi^*}=[0,1,1,1\dots]=\frac{\sqrt{5}-1}{2}$, then in the definition (\ref{defofalpha}) of 
${\mathbb T}_\alpha$ we are always in the first case $(**)$, so that 
${\mathbb T}_{{\Phi^*}}$ is the {\it Fibonacci map} 
$$
{\mathbb T}_{{\Phi^*}}:
[0,1_k,m_{k+1},m_{k+2}, \dots] \in [0,1] \to 
[0,m_{k+1}-1,m_{k+2}, \dots] \in [0,1],
$$
where $1_k$ denotes the string $1,1,\dots 1$ which consists of $k$ consecutive 1's and it is assumed that $n_{k+1}>1$ and $0\leq k <\infty$. 
This map was introduced in \cite{isola2014continued}
and studied in \cite{bonanno2014thermodynamic} from a dynamical point of view. 

In our subsequent  paper \cite{muhammed2015jimm} we also announced some similar and complementary results.
In the next section we show that ${\mathbb T}_{{\Phi^*}}$ and ${\mathbb T}_{{0}}$ are in fact conjugates, albeit via a discontinuous involution of the real line induced by the outer automorphism of the extended modular group $\pgl$, studied in \cite{muhammed2015jimm}.

\begin{figure}[h]
\begin{center}\label{tjimm}
\noindent{\includegraphics[scale=.250]{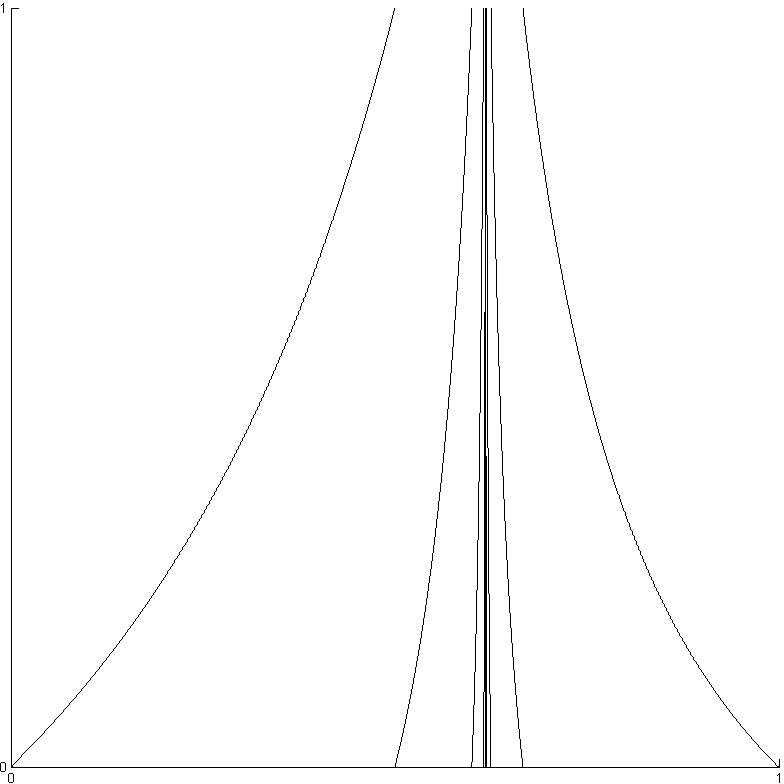}}
\noindent{\includegraphics[scale=.250]{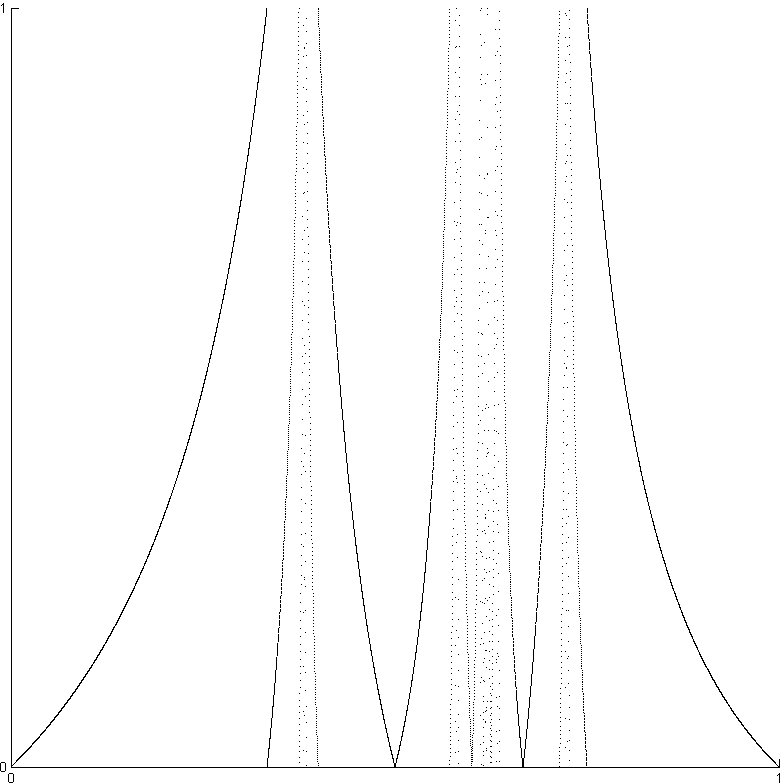}}
\end{center}
\vspace{-.3cm}
\caption{\small The graph of the Fibonacci map ${\mathbb T}_{{\Phi^*}}$ and its first iteration.}
\end{figure}

The map ${\mathbb T}_{{\Phi^*}}$ have the invariant measure $1/x(x+1)$ and gives rise to a transfer operator 
($F_n$ denotes the $n$th Fibonacci number) 
$$
({\mathscr L}_s^{{\Phi^*}}\psi)(y)=\sum_{k=1}^\infty \frac{1}{(F_{k+1}y+F_{k})^{2s}}\psi\left(\frac{F_{k}y+F_{k-1}}{F_{k+1}y+F_{k}}\right),
$$
analytic eigenfunctions of which satisfies the three-term functional equation
$$
\psi(y)=
\frac{1}{y^{2s}}\psi\left(\frac{y+1}{y}\right)+\frac{1}{\lambda}\frac{1}{(y+1)^{2s}}\psi\left(\frac{y}{y+1}\right)
$$
for the eigenvalue $\lambda$.

 The Gauss map is sometimes described as an acceleration (induction) of the Farey map. One may view ${\mathbb T}_{{\Phi^*}}$ as another acceleration of the Farey map.

\subsection{Conjugating the Gauss map to the Fibonacci map}\label{preliminaries}
The Fibonacci map is the conjugate of the Gauss map under an involution. 
\begin{figure}[h]
\begin{center}\label{jimmplot}
\noindent{\includegraphics[scale=1]{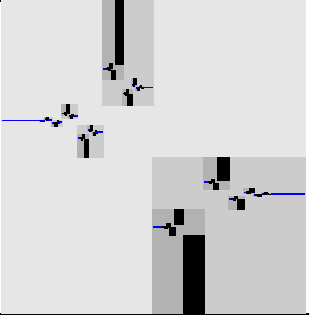}}\\
\caption{\small The plot of Jimm on $[0,1]$.  (The graph lies inside the smaller boxes.)}
\end{center}
\end{figure}
To recall its definition from  \cite{muhammed2015jimm},
 for  $x=[0,n_1,n_2,\dots]$ with  $n_1, \dots \geq 1$, one has
$$
\Jimm(x):=[0,1_{n_1-1},2,1_{n_2-2},2,1_{n_3-2},\dots],
$$
where $1_k$ denotes the sequence ${1,1,\dots, 1}$ of length $k$, and the resulting $1_{-1}$'s are eliminated  in accordance with the rules $[\dots m, 1_{-1},n,\dots]=[\dots m+n-1,\dots]$ and 
$[\dots m, 1_{0},n,\dots]=[\dots m,n,\dots]$. Using the shift description of ${\mathbb T}_0$ and the continued fraction description of the involution $\Jimm$, we see that
$$
{\mathbb T}_{{\Phi^*}}[0,n_1, n_2, \dots]=\Jimm({\mathbb T}_0(\Jimm[0,n_1,n_2,\dots]))
$$
The Fibonacci map has infinitely many fixed points, which are 
$$
{\mathbb T}_{{\Phi^*}}(x)=-\frac{F_{k+1}x-F_{k}}{F_{k+2}x-F_{k+1}}=x\implies x_k=\sqrt{\frac{F_{k}}{F_{k+2}}}=[0,1,\overline{1_{k-1},2}]
$$
with $\lim_{k\rightarrow\infty} x_k=\frac{1}{\Phi}.$
Note that $\Jimm(x_k)=[0,\overline{k+1}]$ with 
$\lim_{k\rightarrow \infty} \Jimm(x_k)=0=\Jimm({\Phi^*})$. The points $\Jimm(x_k)$ are exactly the fixed points of the Gauss map.
More generally, $\Jimm$ gives a correspondence between the periodic orbits of the Gauss map ${\mathbb T}_0$
and those of the map ${\mathbb T}_{{\Phi^*}}$:
$$
{\mathbb T}_0^n(x)=x \iff 
{\mathbb T}_{{\Phi^*}}^n(\Jimm(x))=(\Jimm {\mathbb T}_0\Jimm)^n(\Jimm(x))=(\Jimm {\mathbb T}^n_0\Jimm)(\Jimm(x))=
\Jimm {\mathbb T}^n_0(x)=\Jimm(x)
$$
Hence, $x$ is $n$-periodic orbit for $\mathbb T_0$ if and only if 
$\Jimm(x)$ is $n$-periodic for ${\mathbb T}_{{\Phi^*}}$. 

Finally, note that the following functional equation is satisfied:
$$
{\mathbb T}_{{\Jimm}(\alpha)}(\Jimm x)=\Jimm {\mathbb T}_{{\alpha}}(x).
$$
For $\alpha={\Phi^*}$, this specialises to ${\mathbb T}_{0}(\Jimm x)=\Jimm {\mathbb T}_{{{\Phi^*}}}(x)$ 
$\iff \Jimm{\mathbb T}_{0}(\Jimm x)= {\mathbb T}_{{\Phi^*}}(x)$.

\subsection{Lyapunov exponents}\label{lyapunov}
The {\it Lyapunov exponent} of a map ${\mathbb T}: [0,1] \to [0,1]$
is defined as
\begin{equation} \label{LyapExp} \lambda( x) := \lim_{n\to \infty}
\frac{1}{n} \log | ({\mathbb T}^n)^\prime ( x) | = \lim_{n\to \infty}
\frac{1}{n} \log \prod_{k=0}^{n-1} | {\mathbb T}^\prime ({\mathbb T}^k  x)
|. \end{equation} 
The function $\lambda(x)$ is ${\mathbb T}$--invariant. 
For the Gauss map ${\mathbb T}_0$, the Lyapunov exponent is given by
\begin{equation} \label{LyapT} \lambda_0( x) = 2 \lim_{n\to \infty}
\frac{1}{n} \log q_n( x), \end{equation}
with $q_n( x)$ the successive denominators of the continued
fraction expansion. 
By Khintchin's theorem \cite{khinchin1997continued}
this limit equals
${\pi^2}/{6 \log 2}$ for almost all $x$.
It is known that (see \cite{pollicott1999multifractal}), if the asymptotic proportion of 1's among the partial quotients of $x$ equals 1, 
then the Lyapunov exponent of $x$ (with respect to the Gauss map) is $\lambda_0( x)=2\log \Phi$. 
Now since for almost all $x$, the asymptotic proportion of 1's among the partial quotients of $\Jimm x$  equals 1, we know that 
$\lambda_0(\Jimm x)=2\log \Phi$ for almost all $x$.
It follows that, for almost all $ x$
$$
\lambda_0(\Jimm x)=2\lim_{n\rightarrow\infty} \frac{\log(-C_n\Jimm x+A_n)}{n}=2\log \Phi
$$
$$
\implies \log (-C_n\Jimm x+A_n)  \sim \lambda_0(\Jimm x)=n\log \Phi
$$
$$
\implies -C_n\Jimm x+A_n\sim \exp\left(n\log \Phi)\right)
$$
$$
\implies -C_n x+A_n\sim -C_n x+C_n\Jimm x+\exp\left(n\log \Phi)\right)
$$
$$
\implies 2\lim_{n\rightarrow\infty} \frac{\log(-C_n x+A_n)}{n} \sim 
2\lim_{n\rightarrow\infty}
\log\frac{C_n(\Jimm x- x+\frac{1}{C_n}\exp\left(n\log \Phi\right))}{n}
$$
$$
=2\lim_{n\rightarrow\infty}\frac{\log C_n}{n}+
2\lim_{n\rightarrow\infty}\frac{\log(\Jimm x- x+\frac{1}{C_n}\exp\left(n\log \Phi\right))}{n}
$$
$$
=2\log \Phi+
2\lim_{n\rightarrow\infty}\frac{\log(\Jimm x- x+\frac{1}{C_n}\exp\left(n\log \Phi\right))}{n}
$$
Now there is a bad possibility. The expression inside the brackets is convergent, so it may tend to 0.
However this happens only when 
$$
\Jimm x= x+1\implies \Jimm \Jimm  x=\Jimm( x+1)=1+1/\Jimm x\implies  x=1+\frac{1}{1+ x}\implies  x=\pm \sqrt{2}.
$$
Hence, the Lyapunov exponent of the Fibonacci map is, for almost all $ x$.  
$$
\lambda_{\Phi^*} ( x)=2\log \Phi
$$
again, but for a different reason then the equality $\lambda(\Jimm x)=2\log\Phi$ for almost all $x$. 

The {\it Lyapunov spectrum} of a map is defined as the level sets of the Lyapunov exponent
$\lambda(x)$, by 
$L_c = \{  x \in [0,1] \, | \lambda( x)= c \in \R \}$.
These sets provide a $T$--invariant decomposition of the unit interval,
We don't have a handy description of the Lyapunov spectrum of the map ${\mathbb T}_{Fib}$.

\subsection{Zeta functions}

Since ${\mathscr L}_{s,\alpha=0}(y^t)$ is the Hurwitz zeta $\zeta_{Hur}(y,s+t)$, 
values of the transfer operator ${\mathscr L}_{s, \alpha}$ at the power functions $y^t$ can be viewed as analogues of the Hurwitz zeta:
\begin{eqnarray}\label{dom}
\zeta_{\alpha}(s, t, y):=-\frac{1}{y^{2s+t}}+
\sum_{k=1}^\infty 
\sum_{i=0}^{n_k-1} \left|\frac{\mathrm d}{{\mathrm d}y}
[0,n_1,\dots, n_{k-1},i+y]\right|^{s} [0,n_1,\dots, n_{k-1},i+y]^t
\end{eqnarray}
In particular, 
$$
{\zeta_{\Phi^*}(s, t, y)=
({\mathscr L}_s^{\jimm}\psi)(y^t)=\sum_{k=\mathbf 0}^\infty 
\frac{(F_{k}y+F_{k-1})^t}{(F_{k+1}y+F_{k})^{2s+t}}}
$$
reducing to a two-variable Fibonacci-Hurwitz-zeta for $t=0$
$$
\zeta_{\Phi^*}(s, y):=({\mathscr L}^{s/2}_{\Phi^*} {\mathbf 1})(y)+\frac{1}{y^s}=
\sum_{k=0}^\infty 
\frac{1}{(F_{k+1}y+F_k)^{s}},
$$
and further reducing to the Fibonacci zeta when $y=1$:
$$
\zeta_{\Phi^*}(s):=\sum_{k=1}^\infty 
\frac{1}{F_k^{s}}=\zeta_{\Phi^*}(s, 1)
$$
It is known that the Fibonacci zeta admits a meromorphic continuation to the entire complex plane and its values at negative integers lies in the field $\Q(\sqrt{5})$ (see \cite{murty2013fibonacci} and references therein).
The three-variable zeta satisfies the functional equation 
$$
{\zeta_{\Phi^*}(s, t, 1+1/x)=x^{s-t}\zeta_{\Phi^*}(s, t, x)-x^{-2t}},
$$
which is analogous to the functional equation of the Hurwitz zeta:
 $\zeta_{Hur}(s,1+z)=\zeta_{Hur}(s,z)-(1+z)^{-s}$.

\section{Dynamics of $T_\alpha$ for some special $\alpha$-values.}
Double representation for rationals leads to two different $\mathbb T_\alpha$'s at those points. 
Recall that for $\alpha=q\in \Q$, the two different continued fraction representations of $\alpha$ gives rise to two distinct dynamical maps. 
We shall denote these ${\mathbb T}_{{q^+}}$ when the continued fraction of $q$ ends with $1$ and ${\mathbb T}_{{q^-}}$ otherwise. 

Now, let us discuss some special members of this family of dynamical maps.
\subsection{The map $ {\mathbb T}_1$:} Here we have $\alpha=1=[0,1,\infty]$. 
One has 
$$
\begin{cases} 
x=[0,1,m_{2}, m_{3},m_{4},\dots]&\implies {\mathbb T}_\alpha(x)=[0,m_3,m_4, \dots]={\mathbb T}_0({1\over x}-1)\\
x=[0,m_{1}, m_{2},m_{3},\dots] (m_1>1)& \implies {\mathbb T}_\alpha(x)=[0,m_1-1,m_2, \dots]=({1\over x}-1)^{-1}\\
\end{cases}
$$
\begin{center}
\noindent{\includegraphics[width=14cm, height=4cm]{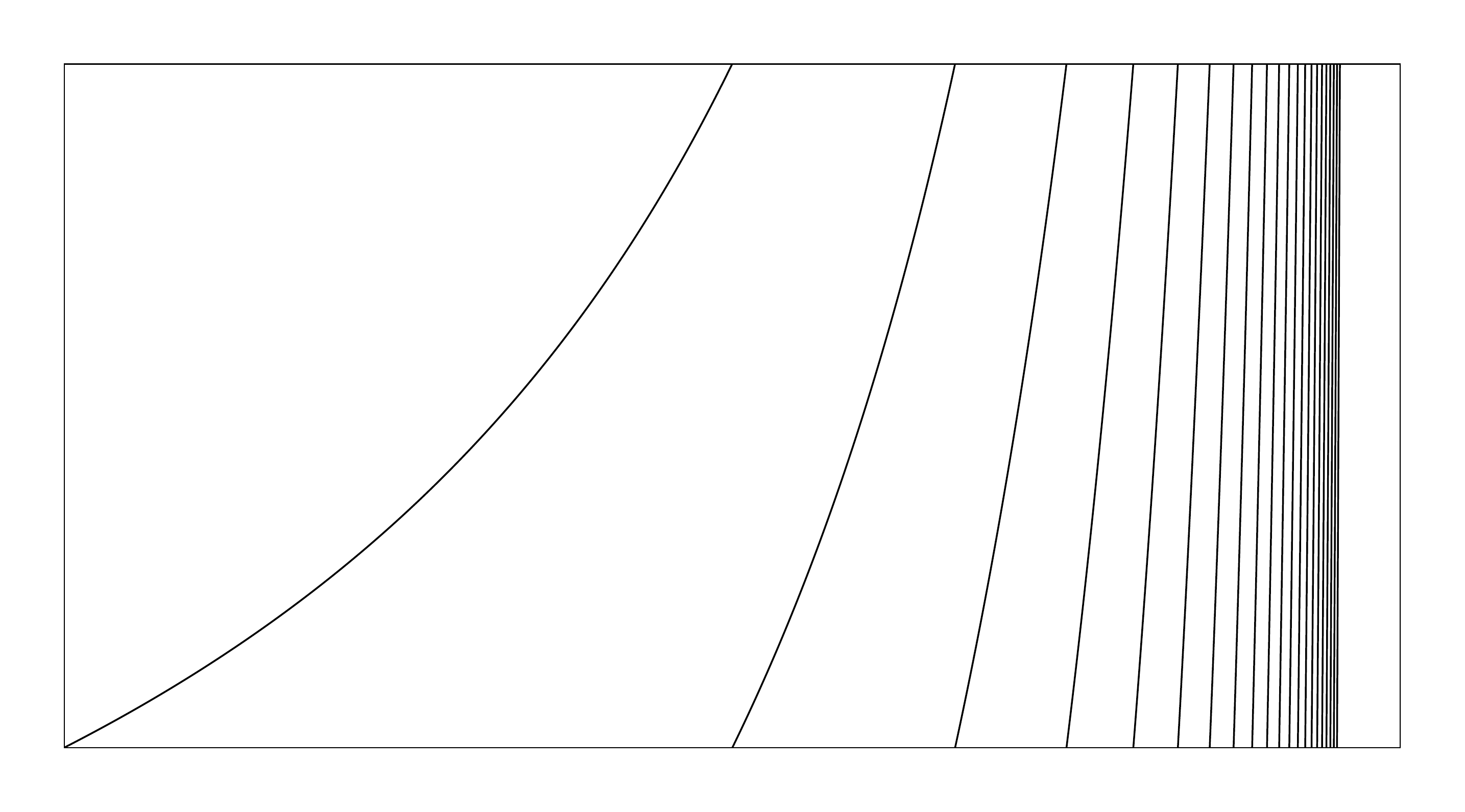}}\\
Graph of the  map $ {\mathbb T}_1$.
\end{center}
One has     
$$
{\mathscr L}_{s, 1} \psi=(\psi|UTU)+\sum_{i=1}^{\infty}(\psi|UTUT^i)
=\sum_{i=0}^{\infty} \frac{1}{(1+i+y)^{2s}}  \psi\left(\frac{i+y}{1+i+y}\right) 
$$
For $s=1$ and $\psi(y)=1/y$ we get 
$$
{\mathscr L}_{s=1, 1} \psi=
\sum_{i=0}^{\infty} \frac{1}{(i+y)(1+i+y)} =
\sum_{i=0}^{\infty} \left(\frac{1}{i+y}- \frac{1}{1+i+y}\right)=\frac{1}{y}
$$
So, $\psi(y)=1/y$ is an invariant measure. 
For the image of the constant function $\psi(y)\equiv 1$ we get the Hurwitz zeta function:
$$
{\mathscr L}_{s, 1} {\mathbf 1}
(s,y)=\sum_{i=0}^{\infty} \frac{1}{(1+i+y)^{2s}}=\zeta_{Hur}(s,y)
$$
We may express the operator ${\mathscr L}_{s, 1}$ in terms of the Mayer transfer operator as
$$
{\mathscr L}_{s, 1} \psi(y)
=\sum_{i=1}^{\infty} \frac{1}{(i+y)^{2s}}  \psi\left(1-\frac{1}{i+y}\right)={\mathscr L}_{s, 0} \widehat{\psi}(y),
$$
where $\widehat{\psi}(y):=\psi(1-1/y)$. In other words, $\psi$ is an eigenfunction for ${\mathscr L}_{s, 1}$ if and only if $\widehat{\psi}$
is an eigenfunction for ${\mathscr L}_{s, 0}$ with the same eigenvalue. 

\subsection{The map $ {\mathbb T}_{1/2}$:} Here we have
$\alpha=1/2^-=[0,2,\infty]$ or $\alpha=1/2^{+}=[0,1,1,\infty]$.
\begin{figure}[h]
\begin{center}\label{t12}
\noindent{\includegraphics[width=14cm, height=4cm]{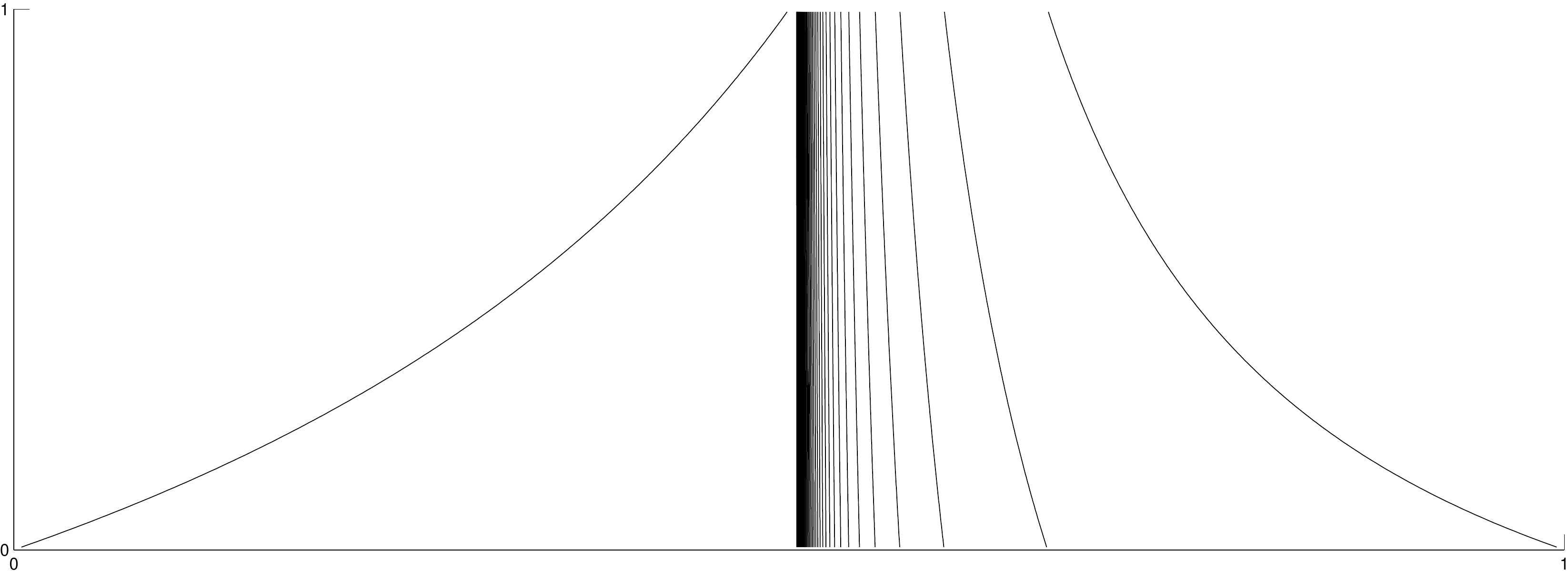}}\\
\caption{\small Plot of the map $ {\mathbb T}_{1/2^+}$}
\end{center}
\end{figure}
\begin{figure}[h]
\begin{center}\label{t12}
\noindent{\includegraphics[width=14cm, height=4cm]{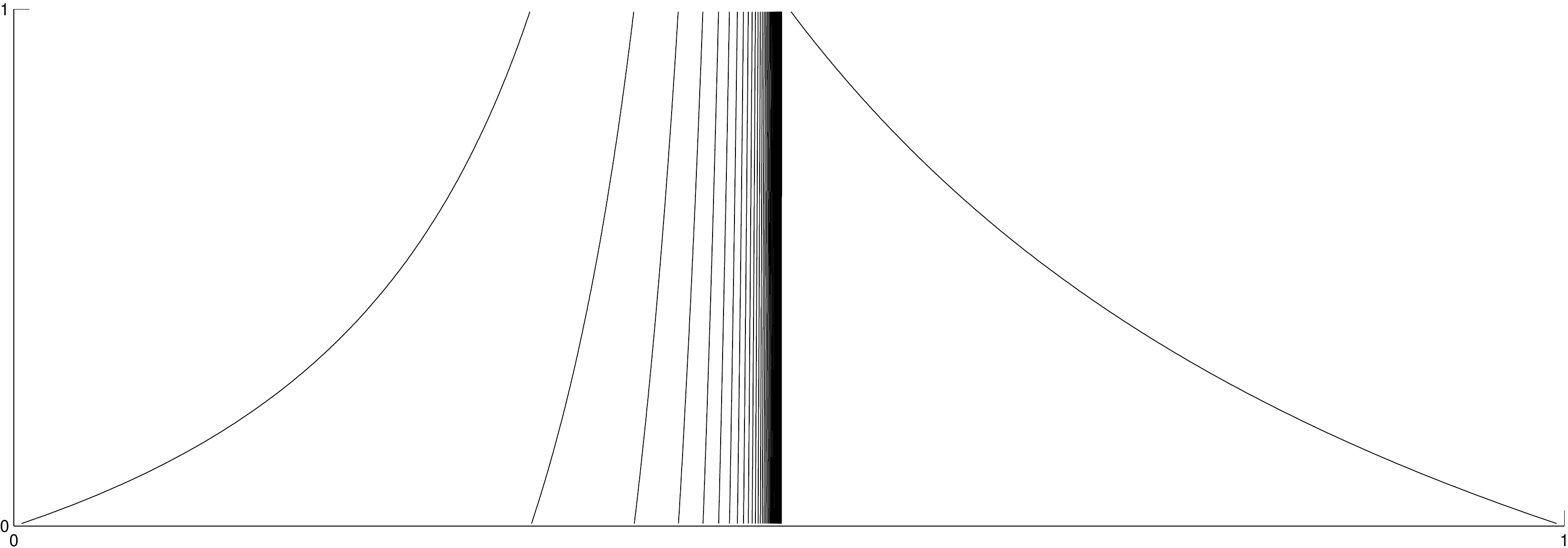}}\\
\caption{\small Plot of the map $ {\mathbb T}_{1/2^-}$}
\end{center}
\end{figure}
  
\begin{eqnarray*}
{\mathscr L}_{s, 1/2^-} \psi=
\frac{1}{(1+y)^{2s}} \psi\left(\frac{1}{1+y}\right)+
\sum_{i=0}^{\infty} 
\frac{1}{(2y+2i+1)^{2s}} \psi\left(\frac{i+y}{2y+2i+1}\right)
\end{eqnarray*}
\begin{eqnarray*}
{\mathscr L}_{s, 1/2^+} \psi=
\frac{1}{(1+y)^{2s}} \psi\left(\frac{y}{1+y}\right)+
\sum_{i=0}^{\infty} 
\frac{1}{(2y+2i+1)^{2s}} \psi\left(\frac{i+y+1}{2y+2i+1}\right)
\end{eqnarray*}
These two operators are equivalent under the transformation 
\begin{eqnarray}\label{equivalencetransf}
{\mathscr L}_{s, 1/2^+} \psi={\mathscr L}_{s, 1/2^-} \overline{\psi}, \mbox{ where } \overline{\psi}(y):=\psi(1-y).
\end{eqnarray}
For the image of the constant function $\psi(y)\equiv 1$ we get:
\begin{eqnarray*}
{\mathscr L}_{s, 1/2^\pm} {\mathbf 1}(s,y)=
\frac{1}{(1+y)^{2s}}+
\sum_{i=0}^{\infty} 
\frac{1}{(2y+2i+1)^{2s}}=\\
\frac{1}{(1+y)^{2s}}+
\sum_{k=1}^{\infty} \frac{1}{(2y+k)^{2s}}
- \sum_{i=1}^{\infty} \frac{1}{(2y+2i)^{2s}}=\\
=\frac{1}{(1+y)^{2s}}+\zeta_{Hur}(2y,2s)-\frac{1}{2^{2s}}\zeta_{Hur}(y,2s)
\end{eqnarray*}
For a fixed function $\psi$ of the operator ${\mathscr L}_{s, 1/2^-}$ one has
\begin{eqnarray*}
\psi(y)=
\frac{1}{(1+y)^{2s}} \psi\left(\frac{1}{1+y}\right)+
\sum_{i=0}^{\infty} 
\frac{1}{(2y+2i+1)^{2s}} \psi\left(\frac{i+y}{2y+2i+1}\right)\implies\\
\psi(y+1)=
\frac{1}{(2+y)^{2s}} \psi\left(\frac{1}{2+y}\right)+
\sum_{i=1}^{\infty} 
\frac{1}{(2y+2i+1)^{2s}} \psi\left(\frac{i+y}{2y+2i+1}\right)
\end{eqnarray*}
Hence, fixed points of the operator ${\mathscr L}_{s, 1/2^-}$ satisfies the functional equation
\begin{eqnarray*}
\psi(y+1)=\psi(y)-\frac{1}{(1+y)^{2s}} \psi\left(\frac{1}{1+y}\right)+
\frac{1}{(2+y)^{2s}} \psi\left(\frac{1}{2+y}\right)-
\frac{1}{(2y+1)^{2s}} \psi\left(\frac{y}{2y+1}\right)
\end{eqnarray*}
In virtue of the equivalence (\ref{equivalencetransf}), 
there is a similar functional equation for fixed functions of the operator ${\mathscr L}_{s, 1/2^+}$. 
Note that we can rewrite the above equation as
\begin{eqnarray*}
\psi(y+1)-\psi(y)+\frac{1}{(1+y)^{2s}}\psi\left(\frac{1}{1+y}\right)=
\frac{1}{(2+y)^{2s}} \psi\left(\frac{1}{2+y}\right)-
\frac{1}{(2y+1)^{2s}} \psi\left(\frac{y}{2y+1}\right)
\end{eqnarray*}
where on the left hand side we have the operator $-\mathscr B_1$.
\subsection{The map $ {\mathbb T}_{1/K}$:} More generally, let us discuss the cases 
$\alpha=1/K^-=[0,K,\infty]$ and $\alpha=1/K^{+}=[0,K-1,1,\infty]$.
\begin{eqnarray*}
{\mathscr L}_{s, 1/K^-} \psi=
\sum_{i=1}^{K-1}\frac{1}{(i+y)^{2s}} \psi\left(\frac{1}{i+y}\right)+
\sum_{i=0}^{\infty} 
\frac{1}{(Ky+Ki+1)^{2s}} \psi\left(\frac{i+y}{Ky+Ki+1}\right)
\end{eqnarray*}
\begin{eqnarray*}
{\mathscr L}_{s, 1/K^+} \psi=
\frac{1}{(y(K-1)+1)^{2s}} \psi\left(\frac{y}{y(K-1)+1}\right)+
\sum_{i=1}^{K-2}\frac{1}{(i+y)^{2s}} \psi\left(\frac{1}{i+y}\right)+\\
\sum_{i=0}^{\infty} 
\frac{1}{(Ky+Ki+K-1)^{2s}} \psi\left(\frac{i+y+1}{Ky+Ki+K-1}\right)
\end{eqnarray*}
These two operators are not visibly equivalent under some transformation.
For a fixed function $\psi$ of the operator ${\mathscr L}_{s, 1/K^-}$ one has
\begin{eqnarray*}
\psi(y)=\sum_{i=1}^{K-1}\frac{1}{(i+y)^{2s}} \psi\left(\frac{1}{i+y}\right)+\sum_{i=0}^{\infty} 
\frac{1}{(Ky+Ki+1)^{2s}} \psi\left(\frac{i+y}{Ky+Ki+1}\right)\implies\\
\psi(y+1)=\sum_{i=2}^{K}\frac{1}{(i+y)^{2s}} \psi\left(\frac{1}{i+y}\right)+\sum_{i=1}^{\infty} 
\frac{1}{(Ky+Ki+1)^{2s}} \psi\left(\frac{i+y}{Ky+Ki+1}\right)
\end{eqnarray*}
Hence, fixed points of the operator ${\mathscr L}_{s, 1/K^-}$ satisfies the functional equation
\begin{eqnarray*}
\psi(y+1)=\psi(y)-\frac{1}{(1+y)^{2s}} \psi\left(\frac{1}{1+y}\right)+
\frac{1}{(K+y)^{2s}} \psi\left(\frac{1}{K+y}\right)-
\frac{1}{(Ky+1)^{2s}} \psi\left(\frac{y}{Ky+1}\right)
\end{eqnarray*}

\subsection{The map $ {\mathbb T}_{\alpha}$ with $\alpha=[0,K,K,\dots]=(-K+\sqrt{K^2+4})/2$:}
\begin{figure}[h]
\begin{center}\label{t12}
\noindent{\includegraphics[width=14cm, height=4cm]{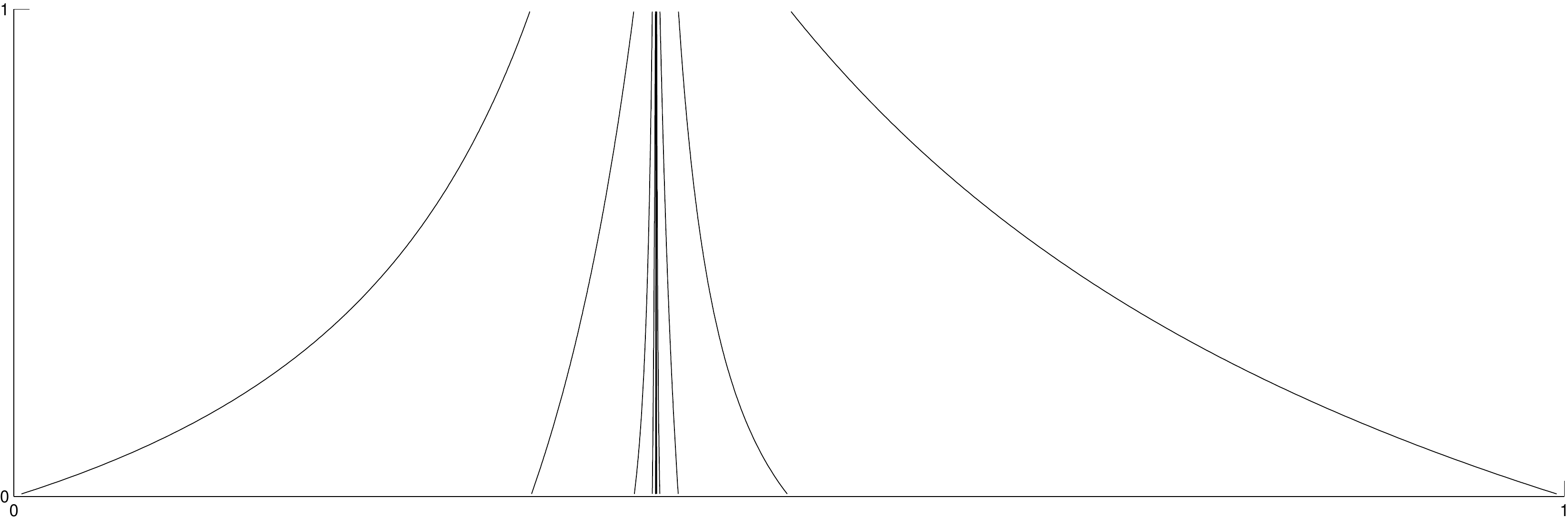}}\\
\caption{\small Plot of the map $ {\mathbb T}_{\sqrt{2}-1}$}
\end{center}
\end{figure}
\begin{figure}[h]
\begin{center}\label{t12}
\noindent{\includegraphics[width=14cm, height=4cm]{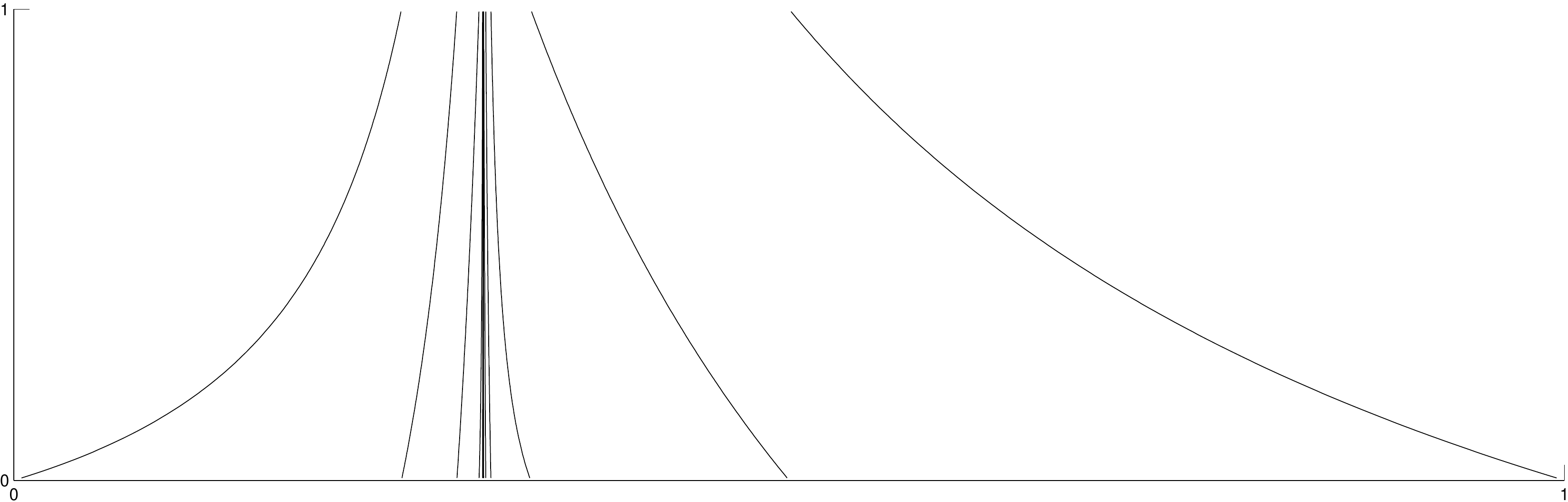}}\\
\caption{\small Plot of the map $ {\mathbb T}_{(\sqrt{13}-3)/2}$}
\end{center}
\end{figure}

\begin{eqnarray*}
{\mathscr L}_{s,\alpha} \psi=
(\psi|-U+
\sum_{i=0}^\infty 
(UT^{K})^iU(I-T)^{-1}(I-T^{K}))\\
=(\psi|-U+
(I-UT^K)^{-1}U(I-T)^{-1}(I-T^{K}))
\end{eqnarray*}
In particular, for $K=1\iff \alpha={\Phi^*}$ we get the transfer operator of the Fibonacci map:
$$
{\mathscr L}_{s, {\Phi^*}}\psi=
(\psi|-U+
\sum_{i=0}^\infty 
(UT)^iU\\
=(\psi|UTU+UTUTU+\dots)
$$

Set
$\varphi:=(\psi|(I-UT^K)^{-1})U(I-T)^{-1}\iff \psi=(\varphi|(I-T)U(I-UT^K)).$
If we assume that $\psi$ is fixed under ${\mathscr L}_{s,\alpha}$, then
\begin{eqnarray*}
\psi=(\psi|-U+(I-UT^K)^{-1}U(I-T)^{-1}(I-T^{K}))\implies\\
(\psi|I+U)=(\psi|(I-UT^K)^{-1}U(I-T)^{-1}(I-T^{K}))\implies\\
(\varphi|(I-T)U(I-UT^K)|(I+U))=(\varphi|I-T^{K})\implies\\
(\varphi|(U-T^K-TU+T^{K+1})(I+U)=(\varphi|I-T^{K})\implies\\
(\varphi|U-T^K-TU+T^{K+1}+I-T^KU-T+T^{K+1}U=(\varphi|I-T^{K})\implies\\
(\varphi|U-TU+T^{K+1}-T^KU-T+T^{K+1}U)=0\implies\\
(\varphi|(T^K-I)(T-U+TU))=0
\end{eqnarray*}
If $\phi$ is $K$-periodic, i.e. $\phi(x+K)=\phi(x)\iff (\phi|T^K)=\phi$, then this equation is satisfied. One has
\begin{align*}
\psi=(\varphi|(I-T)U(I-UT^K))=(\varphi|U-T^K-TU+T^{K+1})=\\
(\varphi|U)-(\varphi|T^K)-(\varphi|TU)+(\varphi|T^{K+1})=\\
(\varphi|U)-(\varphi|I)-(\varphi|TU)+(\varphi|T)=(\varphi|(T-I)(I-U))
\end{align*}
Now, if $\phi$ is not only $K$-periodic, but also $1$-periodic, then this implies that $\psi\equiv 0$. 
Hence we assume that $\phi$ is not 1-periodic $\iff$ $(\phi|T^{K-1}+\dots+T^2+T+I)=0$, i.e.
$$
\phi(x+K-1)+\dots +\phi(x+2)+\phi(x+1)+\phi(x)=0.
$$
Then $\psi:=(\varphi|(T-I)(I-U))$ gives a non-trivial fixed function of the transfer operator.

If $\phi$ is not $K$-periodic, set $\eta:=(\varphi|I-T^K)$.
Then the equation becomes 
\begin{align}\label{kerneleq}
(\eta|T-U+TU)=0 \iff \frac{1}{y^{2s}}\eta\left(\frac{1}{y}\right)=\eta(y+1)+\frac{1}{y^{2s}}\eta\left(1+\frac{1}{y}\right)
\end{align}
which admits $\eta(y)=1/y$ as a solution when $s=1$. Hence, (formally)
$$
\varphi=(\eta|(I-T^K)^{-1})=\sum_{i=0}^\infty (\eta|T^{iK})\implies \varphi(y)=\sum_{i=0}^\infty \frac{1}{y+Ki}.
$$
One has then
\begin{align*}
\psi=(\varphi|U-T^K-TU+T^{K+1})=\\
\sum_{i=0}^\infty 
\frac{1}{y(1+Kiy)}
-\frac{1}{y+Ki+K}
-\frac{1}{y(1+(Ki+1)y)}
+\frac{1}{y+Ki+K+1}\\
=
\sum_{i=0}^\infty \frac{1}{(1+Kiy)(1+Ky+Kiy)}-\frac{1}{(y+Ki+K+1)(y+Ki+K)}
\end{align*}
As a check, for $K=1$ one has, as expected,
\begin{align*}
\psi(y)=
\sum_{i=0}^\infty \frac{1}{y}\left(\frac{1}{1+iy}-\frac{1}{1+(i+1)y)}\right)+\left(\frac{1}{y+i+2}-\frac{1}{y+i+1}\right)\\
= \frac{1}{y}-\frac{1}{y+1}=\frac{1}{y(y+1)}.
\end{align*}
For the case $K=2$ we have
$$
[0,\underbrace{2,2,\dots,2}_{k}, y]=\frac{P_ky+P_{k-1}}{P_{k+1}y+P_{k}},
$$
where $(P_k)_{k=0}^\infty=0, 1, 2, 5, 12, 29, 70, 169, 408, 985, 2378, 5741, 13860, \dots $ 
is the sequence of Pell numbers defined by the recursion $P_0=0$, $P_1=1$ and $P_{k}=2P_{k-1}+P_{k-2}$. 
One has ${\mathscr L}_{s, \alpha} \psi(y)=$
\begin{eqnarray*}
\sum_{i=1}^\infty \frac{1}{(P_{i+1}y+P_{i})^s} \psi\left(\frac{P_iy+P_{i-1}}{P_{i+1}y+P_{i}}\right)+
\sum_{j=1}^\infty 
\frac{1}{(P_{j+1}y+P_{j+1}+P_{j})^s} \psi\left(\frac{P_jy+P_j+P_{j-1}}{P_{j+1}y+P_{j+1}+P_{j}}\right),
\end{eqnarray*}
and the invariant function is
\begin{align*}
\psi(y)=\sum_{i=0}^\infty  \frac{1}{(1+2iy)(1+2y+2iy)}-\frac{1}{(y+2i+3)(y+2i+2)}.
\end{align*}

\bigskip\noindent
{\bf Acknowledgements.} We are grateful to Stefano Isola for commenting on a preliminary version of this paper.
This research is sponsored by the Tübitak grant 115F412 and by a Galatasaray University research grant. 

\begin{center}
\noindent{\includegraphics[scale=.50]{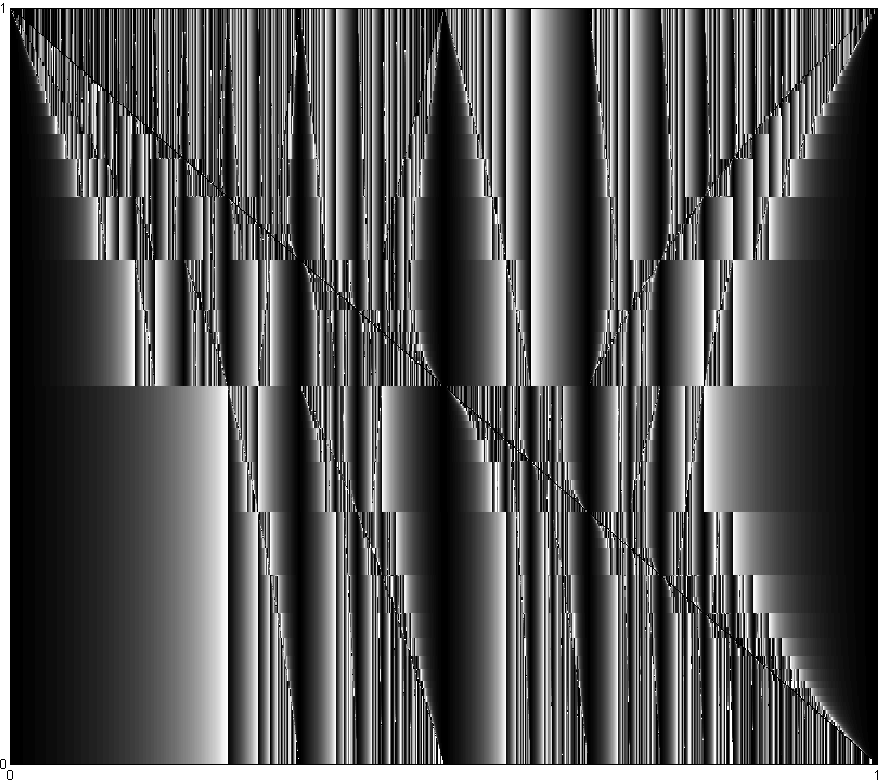}}\\
{\small Third iteration of $ {\mathbb T}_\alpha(x)$.  The intensity is proportional to the value of $T^3_\alpha(x)$.}
\end{center}
%

%\newpage
%REFERENCES TO ADD 
%
%\bibitem{karma}
%Ergodic Theory of Numbers (Carus Mathematical Monographs) 
%by Karma Dajani, Cor Kraaikamp
%
%\bibitem{kramp}
%Metrical Theory of Continued Fractions 
%Marius Iosifescu, Cor Kraaikamp
%
%\bibitem{schweigerpaper}
%Ergodic Theory Of Fibred Systems And Metric Number Theory 
%Fritz Schweiger​
%
%\bibitem{schweigerbook}
%Differentiable equivalence of fractional
%linear maps 
%Fritz Schweiger​
%
%\newpage

%\renewcommand\refname{ \bf Bibliography}
\bibliographystyle{abbrv} 
\bibliography{references4}

\end{document}